\documentclass{article}
\usepackage{graphics}
\usepackage{latexsym}
\usepackage{amsfonts}
\usepackage{amssymb}
\usepackage{amsthm}
\usepackage{amsmath}
\usepackage{hyperref}
\author{Sheng Zhang, Brendan Harding}
\title{Discrete Weierstrass Fourier Transform\\and Experiments}
\date{}
\begin{document}
\newtheorem{theorem}{Theorem}[section]
\newtheorem{proposition}[theorem]{Proposition}
\newtheorem{corollary}[theorem]{Corollary}
\newtheorem{lemma}[theorem]{Lemma}
\maketitle

\begin{abstract}
We established a new method called Discrete Weierstrass Fourier Transform, a faster and more generalized Discrete Fourier Transform, to approximate discrete data. The theory of this method as well as some experiments are analyzed in this paper. In some examples, this method has a faster convergent speed than Discrete Fourier Transform.
\end{abstract}

\section{Introduction}
Barnsley, Harding, Vince, and Viswanathan introduced a notion of Weierstrass Fourier Series to approximate rough functions (see \cite{1}).\\

Key idea of constructing Weierstrass Fourier Series is:
\\1. Construct a linear operator on $L^p(\mathbb{R})$;
\\2. Using this linear operator, transform the classical Fourier basis to a new basis;
\\3. By Gram-Schmidt process, get an orthonormal basis;
\\4. Using this orthonormal basis, do what we can do in classical Fourier Analysis.\\

In this paper, we are going to extend their theory by the following steps:
\\1. Following the key idea in their paper, establish a new expression of Weierstrass Fourier Series, which is more suitable for subsequent discussions in this paper.
\\2. Deduce the notion of Discrete Weierstrass Fourier Transform.
\\3. Provide some numerical examples to test this transform.

\section{Weierstrass Fourier Series}
In classical Fourier Analysis, $\{e^{2\pi ikx}\}_{k\in \mathbb{Z}}$ is a complete orthonormal basis for $L^2 ([0,1])$, and we can approximate functions with Fourier Series. Barnsley, Harding, Vince, and Viswanathan constructed another orthonormal basis for $L^2 ([0,1])$, and established a theory for approximating functions with Weierstrass Fourier Series. A key result in their paper is listed below (see Theorem \ref{dwft 1}).

\begin{theorem}[Barnsley et al.) (see \cite{1}, Corollary 3.1]\label{dwft 1}
Assume that $p \in [1, \infty]$, $a, b\in \mathbb{R}$, $b\ne 0$, and $|a|\ne |b|^{\frac{1}{p}}$. For any $g\in L^p (\mathbb{R})$, there is a unique solution $f\in L^p (\mathbb{R})$ to the equation
\begin{equation}\label{dwft 1 eq 1}
f(x)-af(bx)=g(x),
\end{equation}
and the solution is given by the following series, which are absolutely convergent in $L^p(\mathbb{R})$:
$$f(x)=\left\{\begin{array}{ll}
\sum_{m=0}^{\infty}\,a^{m}\,g(b^{m}x) &\textrm{if $|a|<|b|^{\frac{1}{p}}$}\\
\\
-\sum_{m=1}^{\infty}\,\left (\frac{1}{a})^{m}\,g(\frac{x}{b^m} \right ) & \textrm{if $|a|>|b|^{\frac{1}{p}}$}.
\end{array}\right.$$
\end{theorem}

The following corollary is a direct conclusion from this theorem.
\begin{corollary}\label{dwft 2}
Assume that $a \in [0,1)$ and $b=2$. If $g\in L^\infty (\mathbb{R})$, then (\ref{dwft 1 eq 1}) has a unique solution $f\in L^\infty (\mathbb{R})$, and the solution is given by the following series, which are absolutely convergent in $L^\infty(\mathbb{R})$:
\begin{equation}
f(x)=\sum_{m=0}^{\infty}\,a^{m}\,g(2^{m}x).
\end{equation}
Furthermore, $f\in L^2 ([0,1])$.
\end{corollary}
\begin{proof}
Let $p=\infty$ in Theorem \ref{dwft 1}, and the first part of the corollary is proved. As $f\in L^\infty (\mathbb{R})$, $f\in L^\infty ([0,1])$. So, $f\in L^2 ([0,1])$.
\end{proof}

In paper \cite{1}, Barnsley et al. introduced Weierstrass Fourier Series by substituting $g(x)=\sin kx$ and $\cos kx$. However, another form shown below is more suitable for subsequent discussions in this paper.

In (\ref{dwft 1 eq 1}), assume that $a \in [0,1)$ and $b=2$. For each $k\in \mathbb{Z}$, let $g(x)=e_k(x)=e^{2\pi ikx}$. By Corollary \ref{dwft 2}, there exists a unique solution $f_k(x)=\sum_{m=0}^{\infty}\,a^{m}\,e_k(2^{m}x)\in L^2 ([0,1])$. By normalizing $\{f_k\}$, a normalized basis $\{\hat{e}_k\}$ for $L^2 ([0,1])$ is established:
$$\hat{e}_k=\left\{\begin{array}{ll}
1 &\textrm{if $k=0$}\\
\\
\sqrt{1-a^2}\sum_{m=0}^{\infty}a^m e^{2\pi ik\cdot 2^m x} &\textrm{if $k\ne 0$}.
\end{array}\right.$$

By Gram-Schmidt process, an orthonormal basis $\{\tilde{e}_k\}$ for $L^2 ([0,1])$ is obtained:
$$\tilde{e}_k=
\left\{\begin{array}{lll}
1 &&\textrm{if $k=0$}\\
\\
\hat{e}_k &=\sqrt{1-a^2}\sum_{m=0}^{\infty}a^m e^{2\pi ik\cdot 2^m x} &\textrm{if $k$ is odd}\\
\\
\frac{\hat{e}_k-a\hat{e}_{k/2}}{\sqrt{1-a^2}} &=(1-a^2)\sum_{m=0}^{\infty}a^m e^{2\pi ik\cdot 2^m x}-ae^{\pi ikx} &\textrm{if $k$ is even and $k\ne 0$}.
\end{array}\right.$$

Notice that when $a=0$, $\tilde{e}_k=e_k$ for all $k\in \mathbb{Z}$. To sum up:

\begin{theorem}\label{dwft 3}
Given $a\in [0,1)$, the set of functions $\{\tilde{e}_k\}_{k \in \mathbb{Z}}$ is an orthonormal basis for $L^2 ([0,1])$.
\end{theorem}

Weierstrass Fourier coefficients for a function $f$ can be calculated as follows:
$$\tilde{\alpha}_k=\langle f,\tilde{e}_k \rangle =
\left\{\begin{array}{lll}
\alpha_0 &&\textrm{if $k=0$}\\
\\
\langle f,\hat{e}_k \rangle &=\sqrt{1-a^2}\sum_{m=0}^{\infty}a^m \alpha_{k\cdot 2^m} &\textrm{if $k$ is odd}\\
\\
\frac{\langle f,\hat{e}_k \rangle-a\langle f,\hat{e}_{k/2} \rangle}{\sqrt{1-a^2}} &=(1-a^2)\sum_{m=0}^{\infty}a^m \alpha_{k\cdot 2^m}-a\alpha_{k/2} &\textrm{if $k$ is even and $k\ne 0$},
\end{array}\right.$$
where $\alpha_k$ is the $k^{th}$ Fourier coefficient. Notice that when $a=0$, $\tilde{\alpha}_k=\alpha_k$ for all $k\in \mathbb{Z}$. Now, we can approximate any function in $L^2 ([0,1])$ using Weierstrass Fourier Series.

\begin{theorem}\label{dwft 4}
With the notations above, given $a\in [0,1)$, if $f\in L^2 ([0,1])$ has Fourier expansion
$$f(x)=\sum_{k=-\infty}^{\infty}\alpha_k e_k(x),$$
then it has Weierstrass Fourier expansion
$$f(x)=\sum_{k=-\infty}^{\infty}\tilde{\alpha}_k \tilde{e}_k(x).$$
In particular, when $a=0$, two expansions are the same.
\end{theorem}

\begin{proposition}\label{dwft 5}
With the notations above, given $a\in [0,1)$, if $f\in L^2 ([0,1])$ is a real function, then the following statements hold:\\
(1) $\overline{\tilde{\alpha}_k}=\tilde{\alpha}_{-k}$ for all $k\in \mathbb{Z}$;\\
(2) $\overline{\tilde{e}_k}=\tilde{e}_{-k}$ for all $k\in \mathbb{Z}$;\\
(3) $\sum_{k=-n}^n \tilde{\alpha}_k \tilde{e}_k(x)$ is real and $\sum_{k=-n}^n \tilde{\alpha}_k \tilde{e}_k(x)=\alpha_0 +2Re(\sum_{k=1}^n \tilde{\alpha}_k \tilde{e}_k(x))$ for all $n\in \mathbb{Z_+}$.
\end{proposition}

\section{Discrete Weierstrass Fourier Transform}
With the notations above, given $a\in [0,1)$, for any fixed $n\in \mathbb{N}_+$, let $A$ be an $n\times n$ matrix such that
\begin{equation}
A_{ij}=\tilde{e}_j(\frac{i}{n}), {0\le i,j\le n-1}.
\end{equation}

Given $n$ data points $b_0, b_1, b_2, \cdots, b_{n-1}\in \mathbb{C}$, let $b=(b_0, b_1, b_2, \cdots, b_{n-1})^T$. \textbf{Discrete Weierstrass Fourier Transform}(DWFT) is defined as the linear operator on $\mathbb{C}^n$: $dwft(b)=A^{-1}\,b$. \textbf{Inverse Discrete Weierstrass Fourier Transform}(IDWFT) is defined as the linear operator on $\mathbb{C}^n$: $idwft(b)=A\,b$.
\begin{theorem}\label{dwft 6}
Given $a\in [0,1)$, $dwft$ and $idwft$ definded above are inverse linear operators on $\mathbb{C}^n$. In particular, when $a=0$, $dwft$ is classical Discrete Fourier Transform and $idwft$ is classical Inverse Discrete Fourier Transform.
\end{theorem}

Like DFT, DWFT is useful for data compression. Fix $n\in \mathbb{N}_+$. Let $k\in \mathbb{N}_+$ and $1\le k\le \frac{n+1}{2}$ if $n$ is odd, $1\le k\le \frac{n+2}{2}$ if $n$ is even. Assume
\begin{equation}\label{dwft 6 eq 1}
b=(b_0, b_1, b_2, \cdots, b_{n-1})^T
\end{equation}
is a set of data. Do DWFT on $b$ and get
\begin{equation}\label{dwft 6 eq 2}
c=dwft(b)=(c_0, c_1, c_2, \cdots, c_{n-1})^T.
\end{equation}
Change $c_k, c_{k+1}, \cdots, c_{n-k}$ to $0$ and get
\begin{equation}\label{dwft 6 eq 3}
c'(k)=(c_0, c_1, \cdots, c_{k-1}, 0, \cdots, 0, c_{n-k+1}, \cdots, c_{n-1})^T.
\end{equation}
Do IDWFT on $c'(k)$ and get
\begin{equation}\label{dwft 6 eq 4}
b'(k)=idwft(c'(k))=(b_0', b_1', b_2', \cdots, b_{n-1}')^T.
\end{equation}
When $b=(b_0, b_1, b_2, \cdots, b_{n-1})^T$ is real, $b'(k)$ we got above is not necessarily real. However, we always want to get real approximation for real data. So, if the data are real, a last procedure that discarding the imaginary part of $b'(k)$ should be taken, say
\begin{equation}\label{dwft 6 eq 5}
b''(k)=Re(b'(k)).
\end{equation}
Say $b''(k)$ is the approximation of $b$ using $k$ terms using DWFT. Nevertheless, the imaginary part of $b'(k)$ is usually very small compared to its real part.

Now, if $b''(k)$ is very closed to $b$ even if $k$ is very small compared to $n$, we can store $c'(k)$, $a$, and $n$ instead of data $b$. Then we can calculate matrix $A$ as well as $b''(k)$ when we need and use $b''(k)$ as an approximation of $b$. In this way, plenty of space is saved.

If DWFT and IDWFT are replaced by DFT and IDFT, then the approximation using DFT is obtained.

By Theorem \ref{dwft 6}, DWFT is a generalization of DFT. Since we can choose a suitable $a$ for each set of data, DWFT will never be worse than DFT in approximation.

\section{Experiments}
In this section, some data sets are approximated using DFT and DWFT. With the notations in (\ref{dwft 6 eq 1}), (\ref{dwft 6 eq 2}), (\ref{dwft 6 eq 3}), (\ref{dwft 6 eq 4}) and (\ref{dwft 6 eq 5}), given data vector $b$ and $a\in [0,1)$, \textbf{error vector} for $k$ terms is defined as
\begin{equation}\label{error vector}
\vec{E}_k=b''(k)-b,
\end{equation}
and \textbf{error function} for $f$ is defined as the function
\begin{equation}\label{error function}
\mu (k)=||\vec{E}_k||_2,
\end{equation}
which maps $k$ to the 2-norm of $\vec{E}_k$, for $k\in \{1, 2, \dots, \frac{n}{2}+1\}$.  In particular, $\mu (\frac{n}{2}+1)=0$ by Theorem \ref{dwft 6}.

\subsection{Data from real functions}
In this part, the number of data $n=1024$ for all examples, and each data vector $b=(f(0), f(\frac{1}{n}), f(\frac{2}{n}), \dots, f(\frac{n-1}{n}))^T$ for some real-valued function $f$ on $[0,1]$.

\subsubsection{Linear function}
\begin{center}
$f(x)=x-0.5$ and $a=0.5$\\
\scalebox{0.5}{\includegraphics{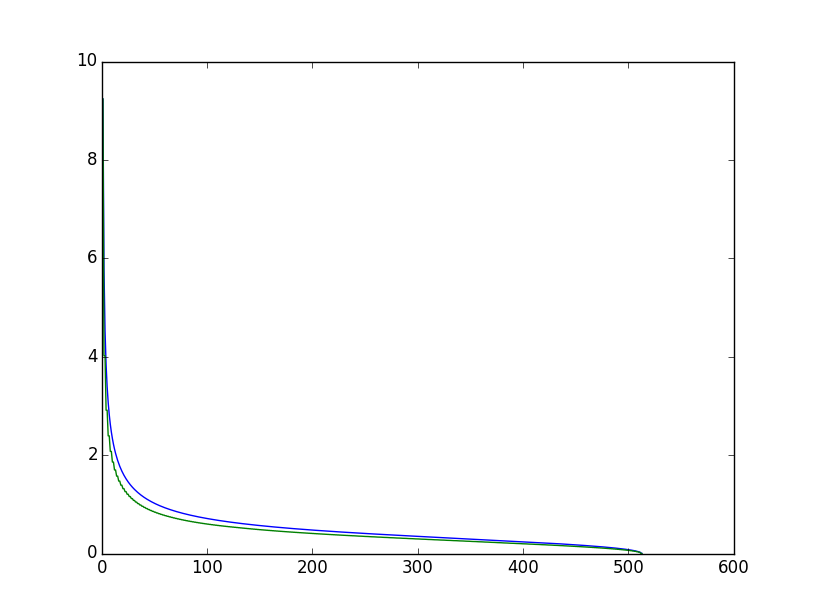}}
\end{center}
Figure 1.1.1  The figure of the error function for $f(x)=x-0.5$ with $a=0.5$.  The blue curve represents approximation using DFT, and the green curve represents approximation using DWFT.  DWFT is always better.

\begin{center}\scalebox{0.5}{\includegraphics{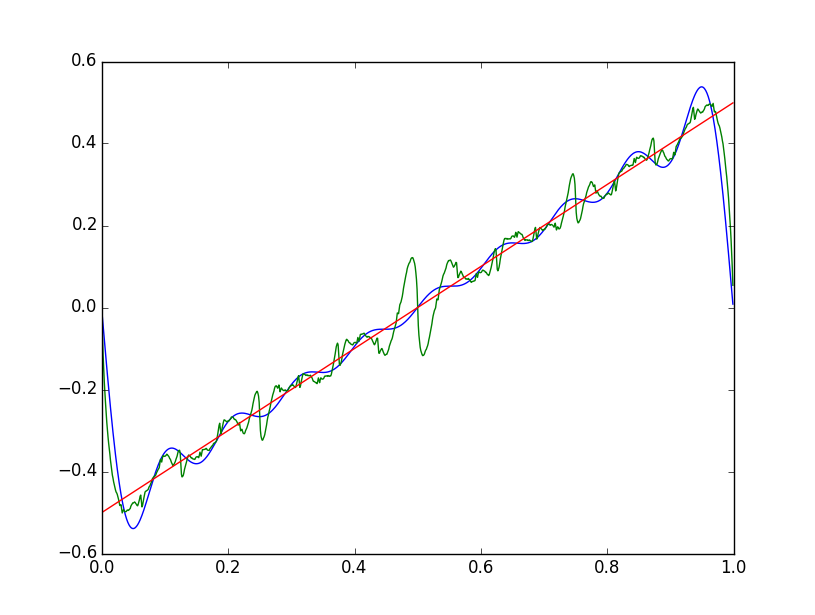}}\end{center}
Figure 1.1.2  Approximation of $f(x)=x-0.5$ using $10$ terms.  The red curve represents the original data, the blue curve represents approximation values using DFT, and the green curve represents approximation values using DWFT.  See also Figure 1.1.3 and Figure 1.1.4.

\begin{center}\scalebox{0.5}{\includegraphics{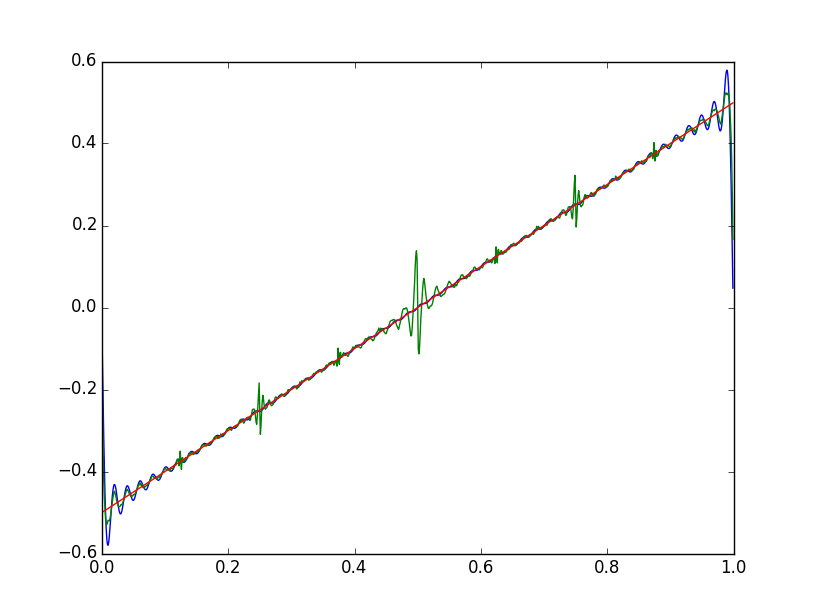}}\end{center}
\begin{center}
Figure 1.1.3  Approximation of $f(x)=x-0.5$ using $50$ terms.
\end{center}

\begin{center}\scalebox{0.5}{\includegraphics{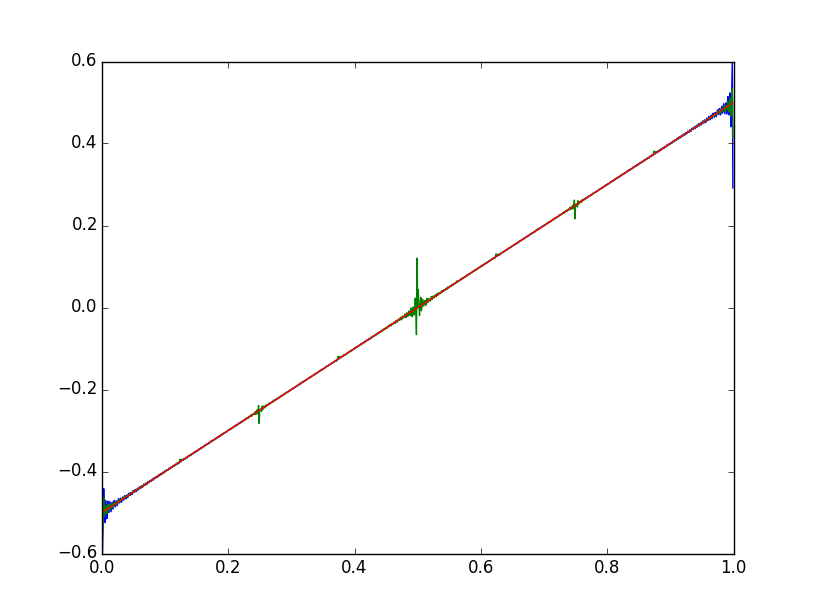}}\end{center}
\begin{center}
Figure 1.1.4  Approximation of $f(x)=x-0.5$ using $300$ terms.
\end{center}

\subsubsection{Triangular function with low frequency and high frequency}
\begin{center}
$f(x)=\sin x+0.01\sin 105x$ and $a=0.5$\\
\scalebox{0.5}{\includegraphics{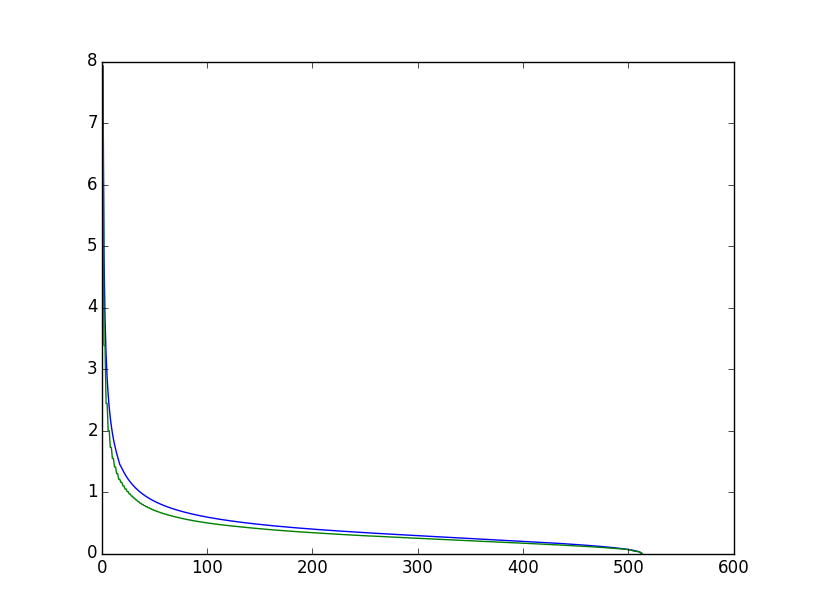}}
\end{center}
Figure 1.2.1  The figure of the error function for $f(x)=\sin x+0.01\sin 105x$ with $a=0.5$.  The blue curve represents approximation using DFT, and the green curve represents approximation using DWFT.  DWFT is better at the first $510$ terms.

\begin{center}\scalebox{0.5}{\includegraphics{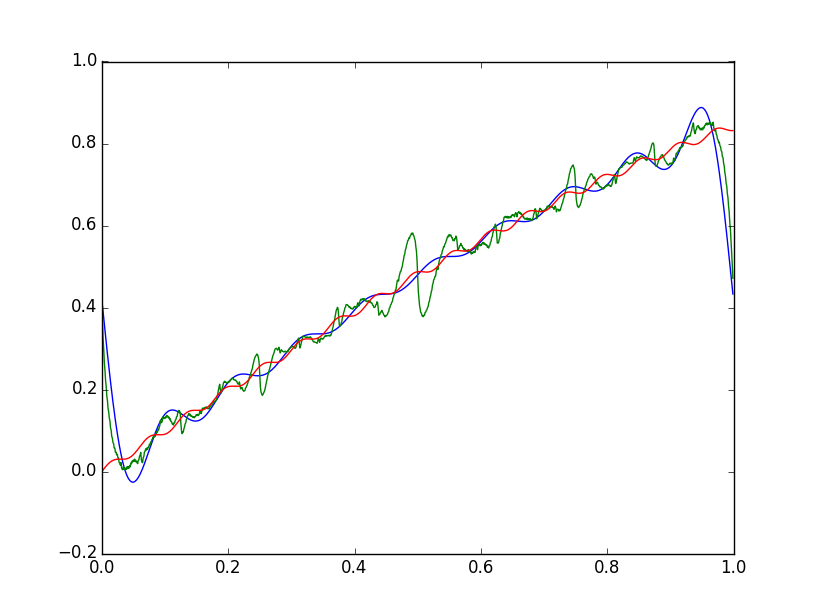}}\end{center}
Figure 1.2.2  Approximation of $f(x)=\sin x+0.01\sin 105x$ using $10$ terms.  The red curve represents the original data, the blue curve represents approximation values using DFT, and the green curve represents approximation values using DWFT.  See also Figure 1.2.3 and Figure 1.2.4.

\begin{center}\scalebox{0.5}{\includegraphics{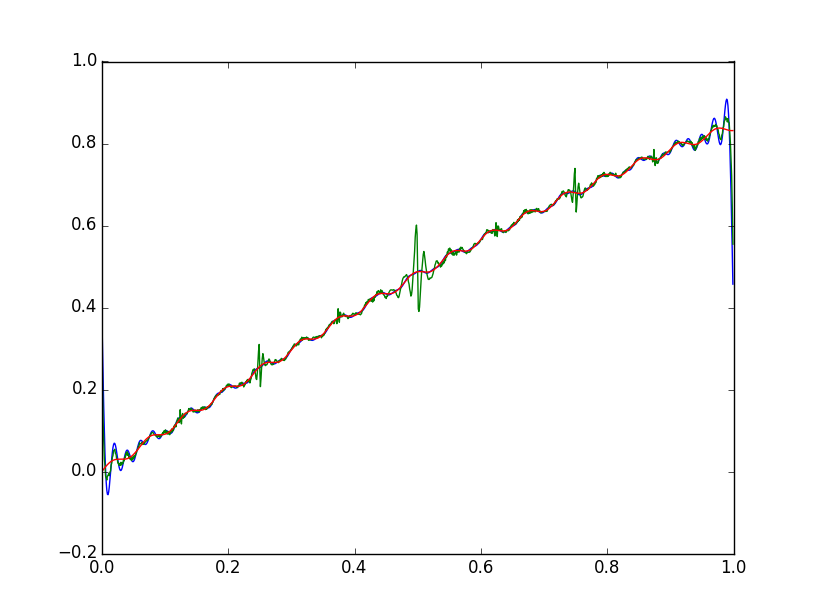}}\end{center}
Figure 1.2.3  Approximation of $f(x)=\sin x+0.01\sin 105x$ using $50$ terms.

\begin{center}\scalebox{0.5}{\includegraphics{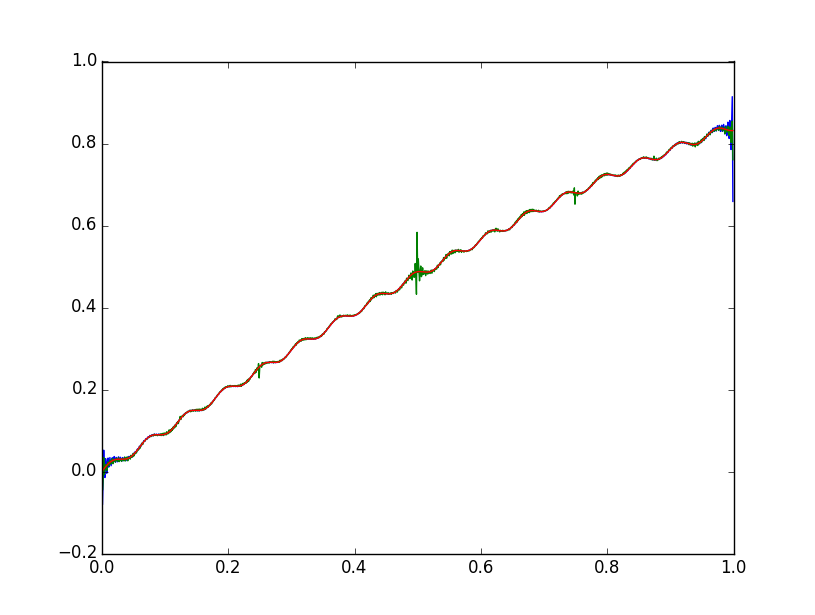}}\end{center}
Figure 1.2.4  Approximation of $f(x)=\sin x+0.01\sin 105x$ using $300$ terms.

\subsubsection{Discontinuous function}
\begin{center}
$f(x)=
\left\{\begin{array}{ll}
0 & x\in [0,\frac{1}{2}]\\
1 & x\in (\frac{1}{2},1]
\end{array}\right.$ and $a=0.5$\\
\scalebox{0.5}{\includegraphics{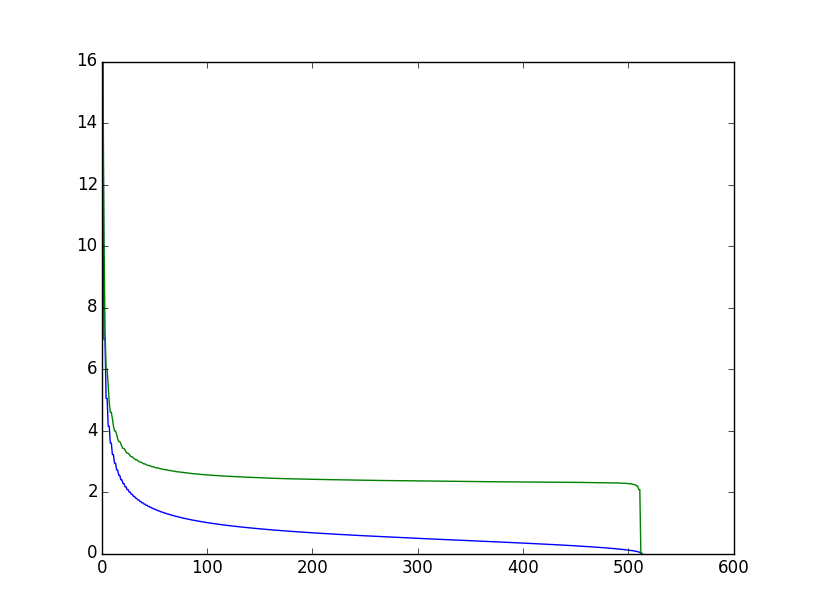}}
\end{center}
Figure 1.3.1  The figure of the error function for $f(x)=
\left\{\begin{array}{ll}
0 & x\in [0,\frac{1}{2}]\\
1 & x\in (\frac{1}{2},1]
\end{array}\right.$
with $a=0.5$.  The blue curve represents approximation using DFT, and the green curve represents approximation using DWFT.  DWFT is not better.

\begin{center}\scalebox{0.5}{\includegraphics{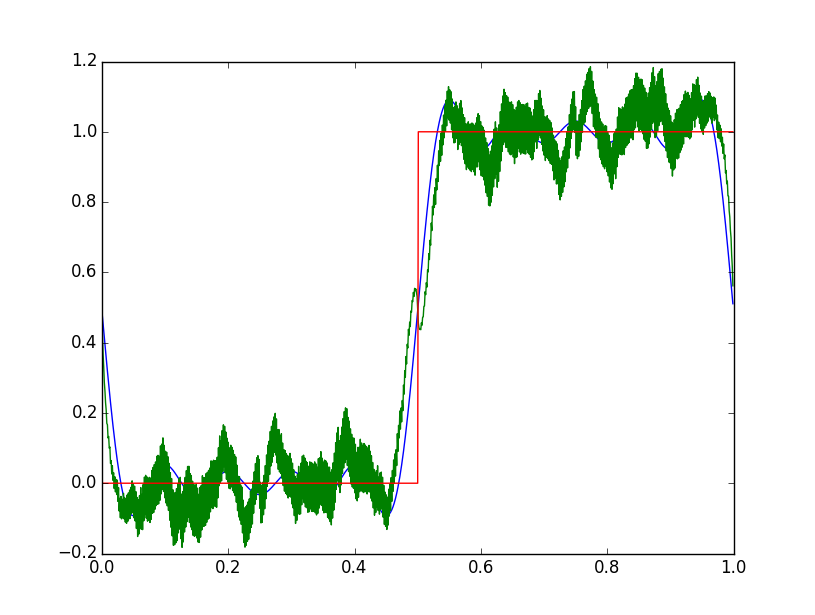}}\end{center}
Figure 1.3.2  Approximation of $f(x)=
\left\{\begin{array}{ll}
0 & x\in [0,\frac{1}{2}]\\
1 & x\in (\frac{1}{2},1]
\end{array}\right.$
using $10$ terms.  The red curve represents the original data, the blue curve represents approximation values using DFT, and the green curve represents approximation values using DWFT.  See also Figure 1.3.3 and Figure 1.3.4.

\begin{center}\scalebox{0.5}{\includegraphics{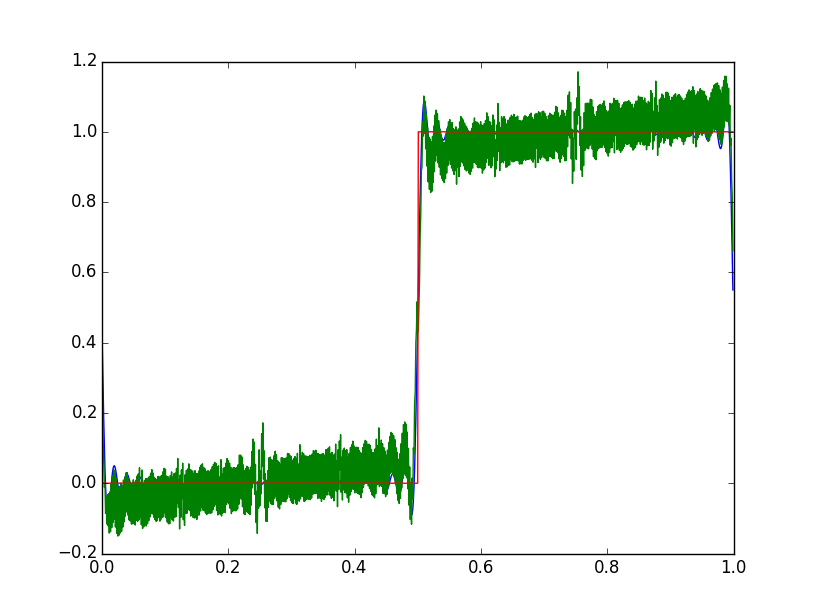}}\end{center}
Figure 1.3.3  Approximation of $f(x)=
\left\{\begin{array}{ll}
0 & x\in [0,\frac{1}{2}]\\
1 & x\in (\frac{1}{2},1]
\end{array}\right.$
using $50$ terms.

\begin{center}\scalebox{0.5}{\includegraphics{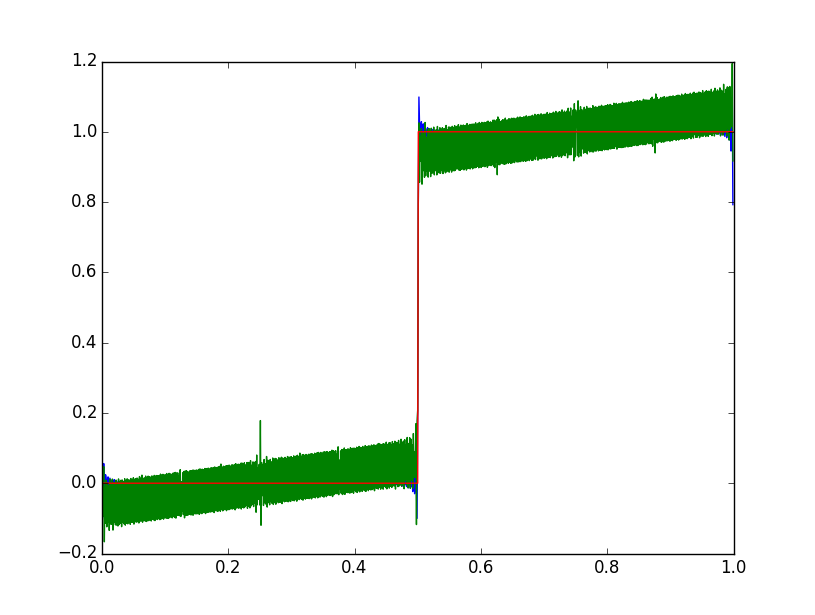}}\end{center}
Figure 1.3.4  Approximation of $f(x)=
\left\{\begin{array}{ll}
0 & x\in [0,\frac{1}{2}]\\
1 & x\in (\frac{1}{2},1]
\end{array}\right.$
using $300$ terms.

\subsubsection{Rough function}
\begin{center}
$f(x)=\sum_{k=0}^{\infty}0.42^k\cos (\pi \cdot 2^k x)$ and $a=0.42$\\
\scalebox{0.5}{\includegraphics{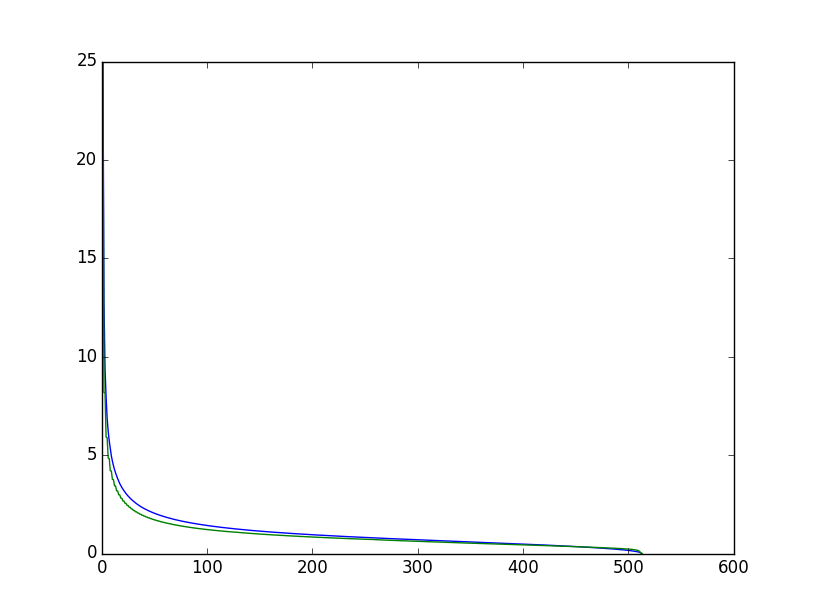}}
\end{center}
Figure 1.4.1  The figure of the error function for $f(x)=\sum_{k=0}^{\infty}0.42^k\cos (\pi \cdot 2^k x)$ with $a=0.42$.  The blue curve represents approximation using DFT, and the green curve represents approximation using DWFT.  DWFT is better at the first $454$ terms.

\begin{center}\scalebox{0.5}{\includegraphics{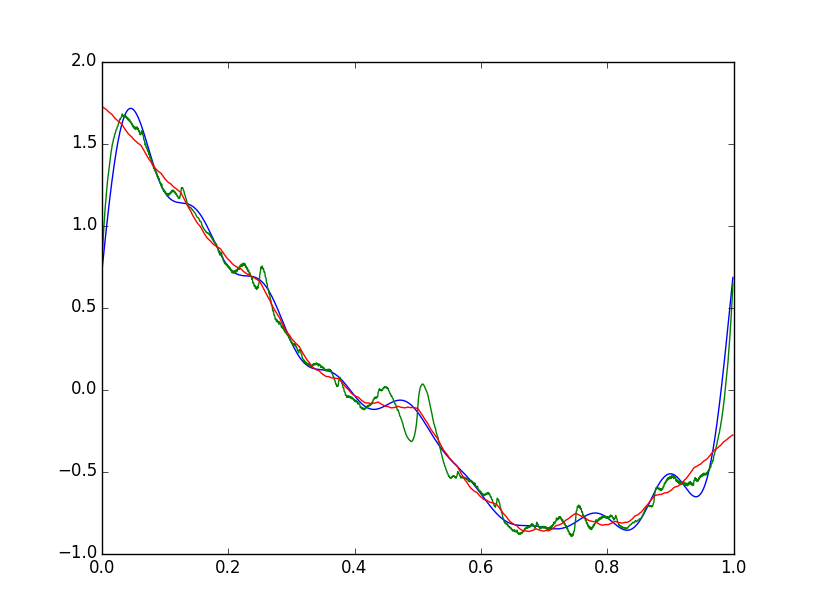}}\end{center}
Figure 1.4.2  Approximation of $f(x)=\sum_{k=0}^{\infty}0.42^k\cos (\pi \cdot 2^k x)$ using $10$ terms.  The red curve represents the original data, the blue curve represents approximation values using DFT, and the green curve represents approximation values using DWFT.  See also Figure 1.4.3 and Figure 1.4.4.

\begin{center}\scalebox{0.5}{\includegraphics{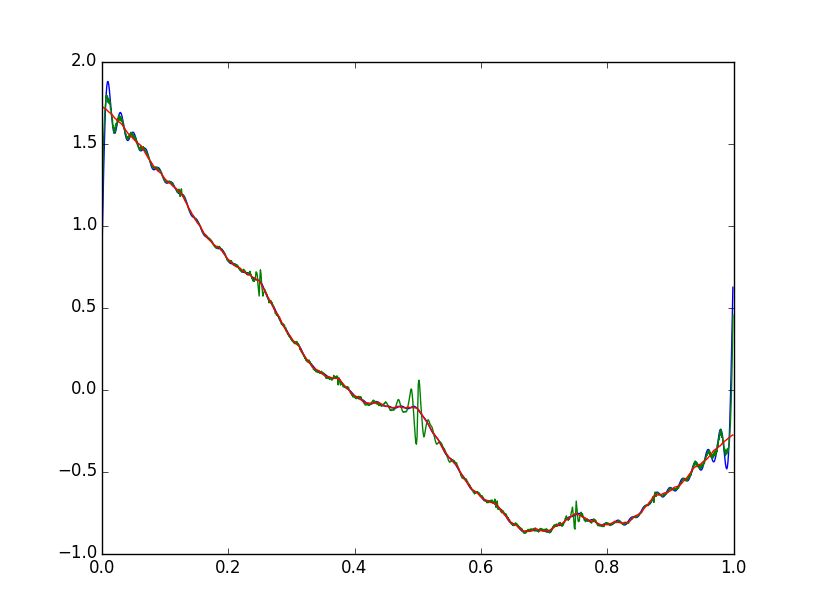}}\end{center}
Figure 1.4.3  Approximation of $f(x)=\sum_{k=0}^{\infty}0.42^k\cos (\pi \cdot 2^k x)$ using $50$ terms.

\begin{center}\scalebox{0.5}{\includegraphics{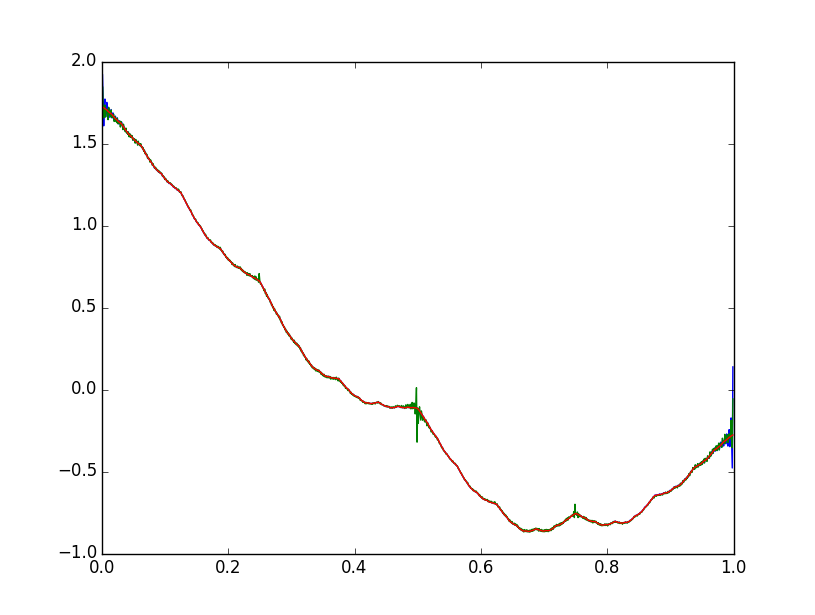}}\end{center}
Figure 1.4.4  Approximation of $f(x)=\sum_{k=0}^{\infty}0.42^k\cos (\pi \cdot 2^k x)$ using $300$ terms.

\subsection{Discrete data}
In this part, data are practical data from some websites. Given data vector $b$ and $a\in [0,1)$, error vector and error function for the data are defined in (\ref{error vector}) and (\ref{error function}).

\subsubsection{Stock price of Commonwealth Bank of Australia}
Daily open prices of Commonwealth Bank of Australia from Sep 30, 2010.  Totally 1024 data points.  Source from https://au.finance.yahoo.com/.
\begin{center}\scalebox{0.5}{\includegraphics{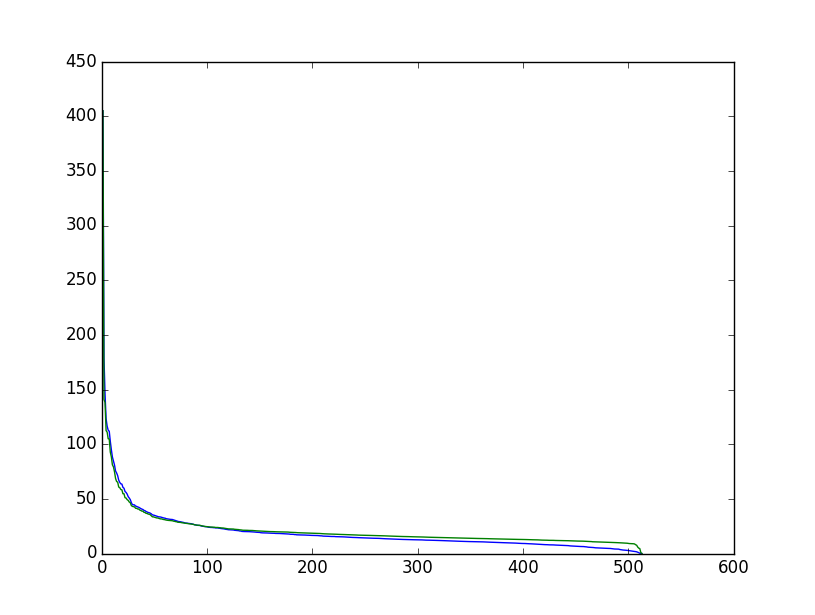}}\end{center}
Figure 2.1.1  The figure of the error function for 1024 daily open prices of Commonwealth Bank of Australia with $a=0.3$.  The blue curve represents approximation using DFT, and the green curve represents approximation using DWFT.  DWFT is better at the first $88$ terms.

\begin{center}\scalebox{0.5}{\includegraphics{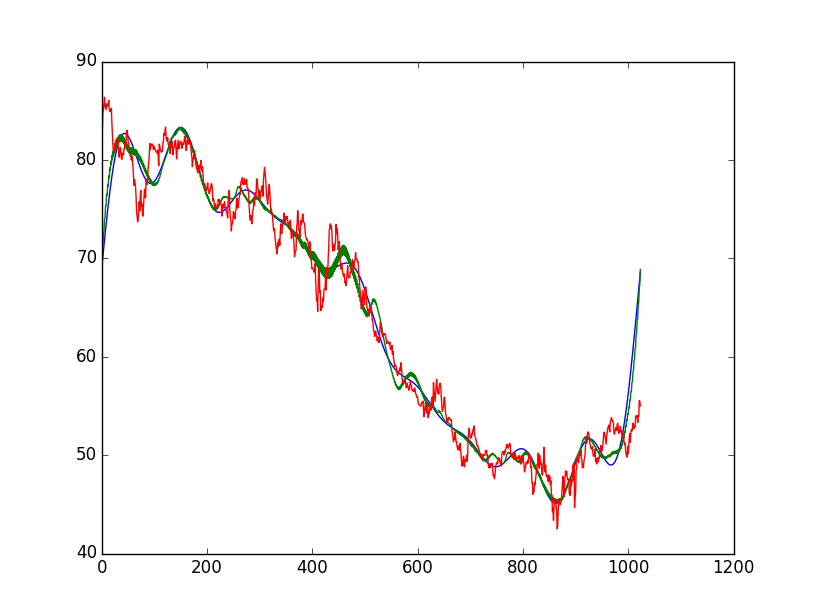}}\end{center}
Figure 2.1.2  Approximation of 1024 daily open prices of Commonwealth Bank of Australia using $10$ terms.  The red curve represents the original data, the blue curve represents approximation values using DFT, and the green curve represents approximation values using DWFT.  See also Figure 2.1.3 and Figure 2.1.4.

\begin{center}\scalebox{0.5}{\includegraphics{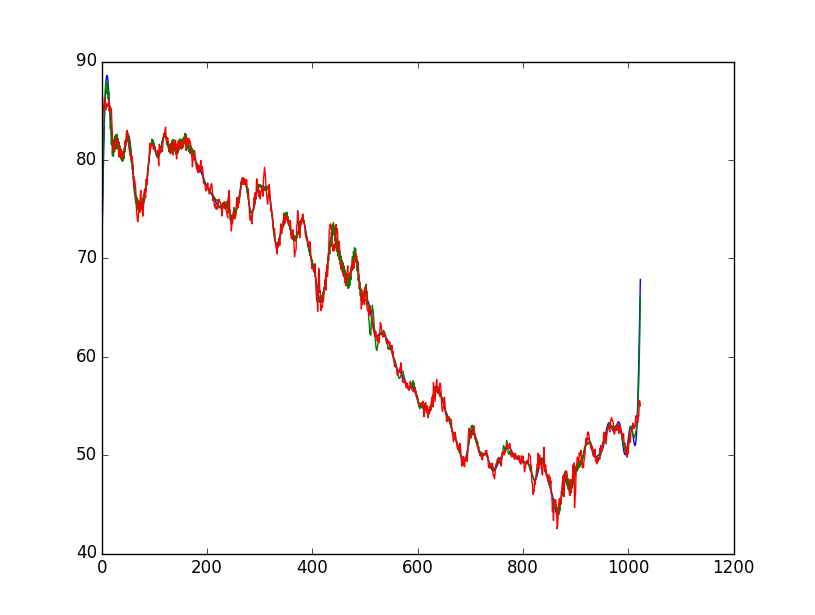}}\end{center}
Figure 2.1.3  Approximation of 1024 daily open prices of Commonwealth Bank of Australia using $50$ terms.

\begin{center}\scalebox{0.5}{\includegraphics{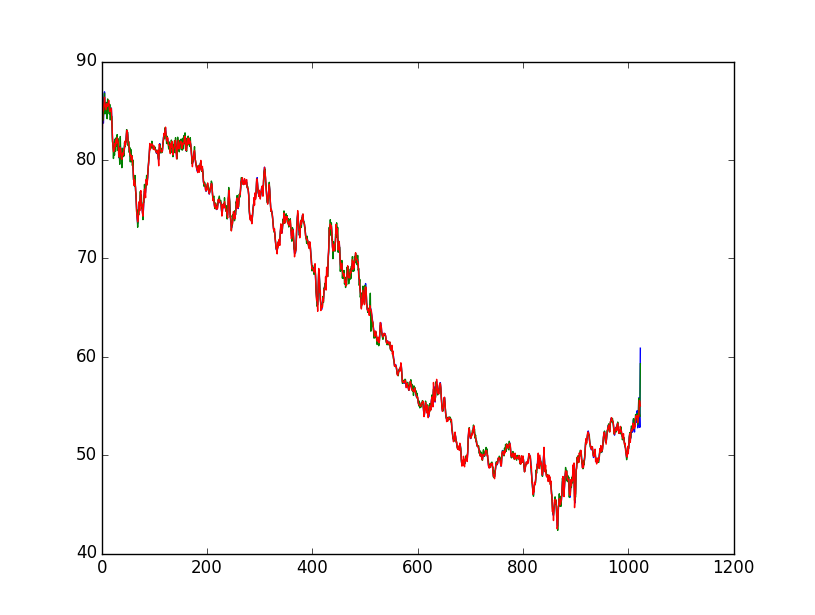}}\end{center}
Figure 2.1.4  Approximation of 1024 daily open prices of Commonwealth Bank of Australia using $300$ terms.

\subsubsection{Water level of Alameda in California}
Water level of Alameda in California per hour from May 29, 2014.  Totally 1024 data points.  Source from http://tidesandcurrents.noaa.gov/.
\begin{center}\scalebox{0.5}{\includegraphics{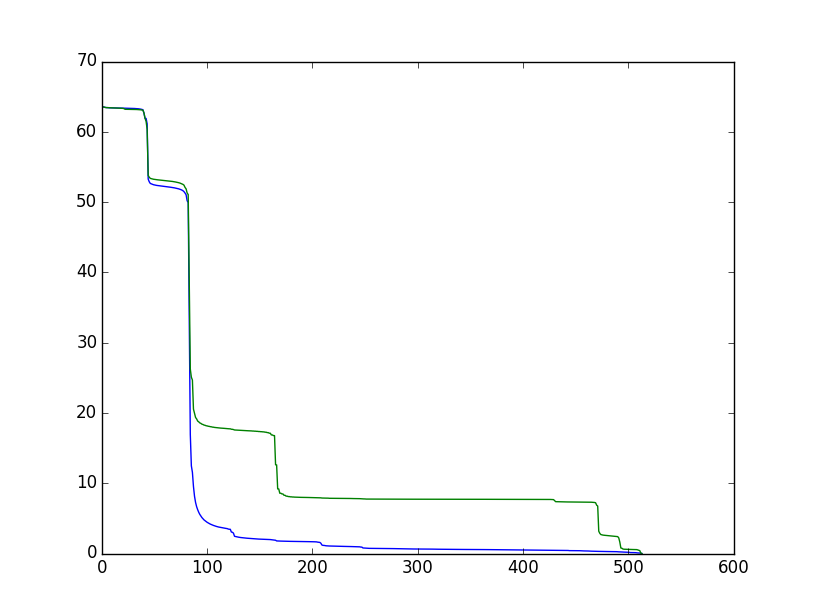}}\end{center}
Figure 2.2.1  The figure of the error function for 1024 water levels of Alameda in California per hour from May 29, 2014 with $a=0.3$.  The blue curve represents approximation using DFT, and the green curve represents approximation using DWFT.  DWFT is not better.

\begin{center}\scalebox{0.5}{\includegraphics{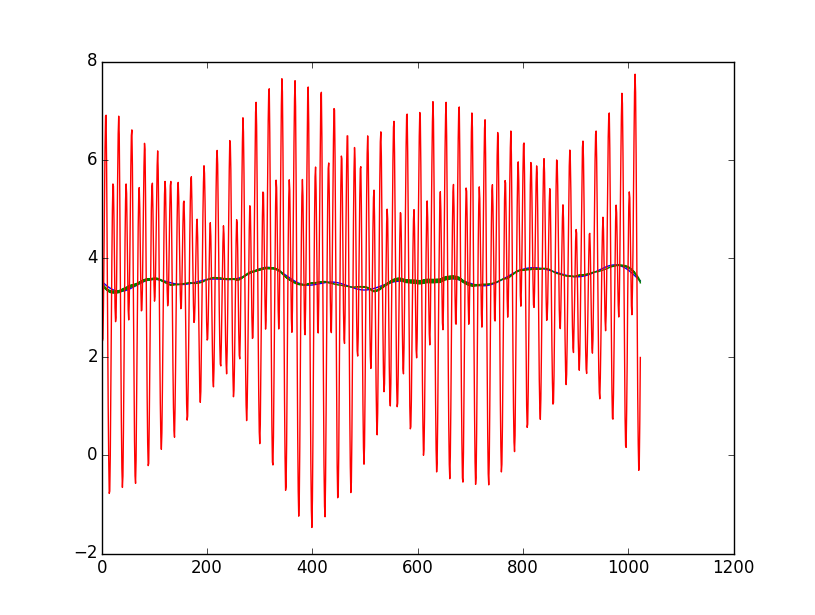}}\end{center}
Figure 2.2.2  Approximation of 1024 water levels of Alameda in California using $10$ terms.  The red curve represents the original data, the blue curve represents approximation values using DFT, and the green curve represents approximation values using DWFT.  See also Figure 2.2.3 and Figure 2.2.4.

\begin{center}\scalebox{0.5}{\includegraphics{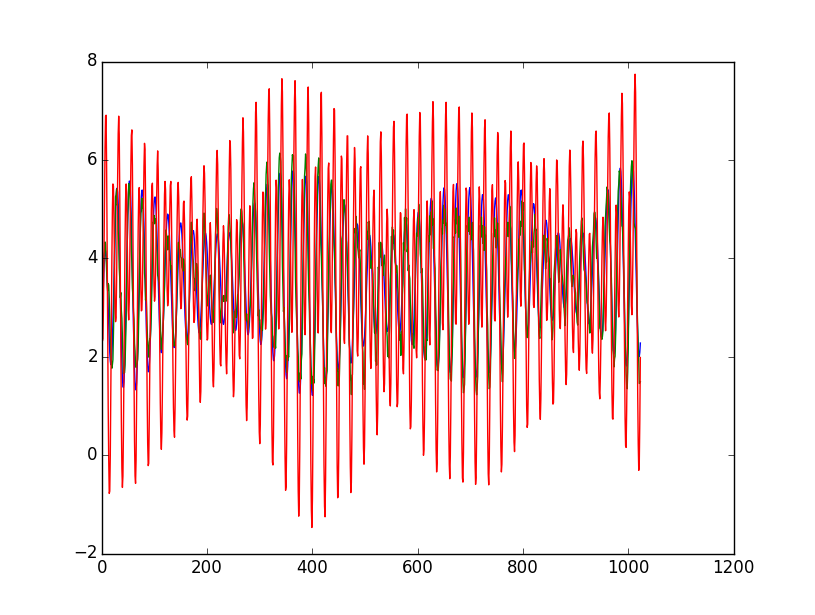}}\end{center}
Figure 2.2.3  Approximation of 1024 water levels of Alameda in California using $50$ terms.

\begin{center}\scalebox{0.5}{\includegraphics{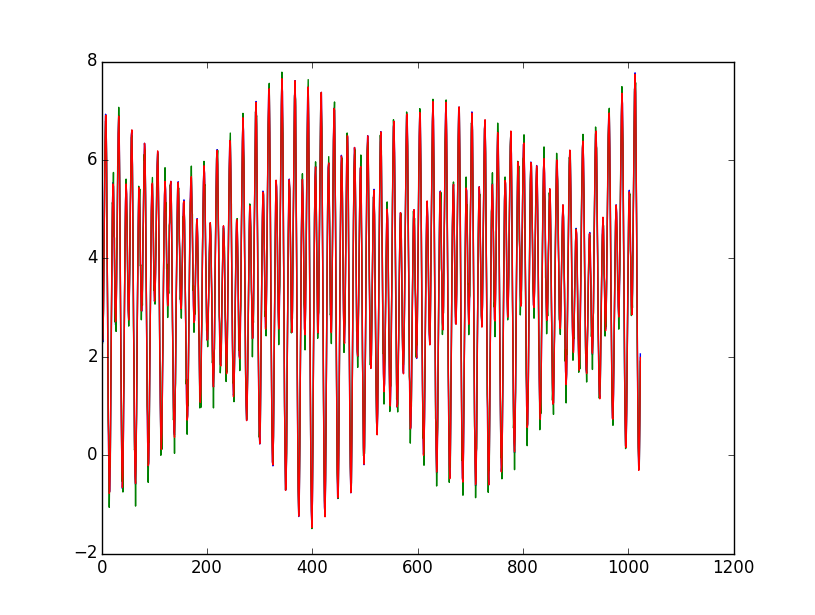}}\end{center}
Figure 2.2.4  Approximation of 1024 water levels of Alameda in California using $300$ terms.

\subsection{Comments}
\subsubsection{}
Figure 1.1.2 and Figure 1.1.3 in Section 4.1.1 coincide with Figure 1, Figure 2, Figure 3 and Figure 4 of \cite{1}, which means that the discrete algorithm in this paper gives the same result as the continuous algorithm when the number of the data is much bigger than the terms of approximation.

\subsubsection{}
In the figures above we can see that at some of the data points, DWFT approximates better, while at some other data points DFT is better.  For example, in Figure 1.1.3, at the data points near 0.0 or 1.0, DWFT is better while at the data points near 0.5, it is significantly worse, but overall, it approximates the data better in the norm defined above.

\subsubsection{}
We can also see a characteristic that approximation values using DWFT shakes heavily at the points near $\frac{1}{2}$, less heavily at the points near $\frac{1}{4},\frac{3}{4}$, and less heavily at the points near $\frac{1}{8},\frac{3}{8},\frac{5}{8},\frac{7}{8}$, $\cdots$

\subsubsection{}
In the above examples, DWFT is better than DFT for some smooth continuous functions.  For some other smooth functions such as $f(x)=x^2$ with $a=0.5$ and $f(x)=e^x$ with $a=0.5$, DWFT is also better than DFT.  But for $f(x)=x(x-1)$ with $a=0.5$, DWFT behaves very bad.

\subsubsection{}
For triangle functions with low frequency and high frequency such as the example in Section 4.1.2, $f(x)=\sin x+0.01\cos 100x+0.01\cos 200x$ with $a=0.5$ and $f(x)=10\sin 0.1x-\cos x+0.01\sin 100x+0.02\cos 200x$ with $a=0.5$, DWFT is better than DFT.

However, make sure that the low frequencies are low enough with coefficients big enough compared to the high frequencies, and the low frequencies are not counteracted each other.  For instance, DWFT doesn't work better than DFT for the function $f(x)=10\sin 0.1x-\sin x+0.01\sin 100x+0.02\cos 200x$ with $a=0.5$ because the low frequency terms ``$10\sin 0.1x$'' and ``$-\sin x$'' are counteracted each other.

\subsubsection{}
For discontinuous functions, DWFT doesn't seem better than DFT.  As we can see in Section 4.1.3, approximation values using DWFT for discontinuous function shake heavily.

\subsubsection{}
As we can see in Section 4.1.4, DWFT has some advantages in approximating self-similar rough functions.

\subsubsection{}
As we can see in Section 4.2, DWFT seems not good for approximating practical data, because practical data are not perfectly self-similar.  A small shake in the data can induce shakes everywhere in the approximation values using DWFT.  In this case, the approximation values using DWFT look much rougher than the data.

\subsubsection{}
If some data set can be approximated better using DWFT, we have a better way to compress the data.  However, the calculation of DWFT is much more complicated than DFT.  We might wish to find a fast way to calculate DWFT.

\end{document}